\documentclass[reqno]{amsart}
\usepackage{amsmath,amsthm,amscd,amssymb,amsfonts, amsbsy}
\usepackage{latexsym, color, enumerate}
\usepackage{pxfonts}
\usepackage{mathrsfs}
\usepackage{marginnote}
\usepackage{todonotes}
\usepackage{hyperref}

\usepackage{mathtools}
\mathtoolsset{showonlyrefs}

\usepackage{thmtools}
\usepackage{thm-restate}

\usepackage{tikz}
\usetikzlibrary{snakes}
\usepackage{multicol}

\usepackage{cancel}

\numberwithin{equation}{section}

\theoremstyle{plain}
\newtheorem{theorem}{Theorem}[section]
\newtheorem{lemma}[theorem]{Lemma}

\theoremstyle{definition}

\theoremstyle{remark}
\newtheorem{remark}[theorem]{Remark}

\newtheorem{conjecture}[theorem]{Conjecture}
\newtheorem{problem}[theorem]{Problem}

\newcommand{\bN}{\mathbb{N}}

\newcommand{\bR}{\mathbb{R}}

\renewcommand{\vec}[1]{\boldsymbol{#1}}

\def\XXint#1#2#3{{\setbox0=\hbox{$#1{#2#3}{\int}$}
		\vcenter{\hbox{$#2#3$}}\kern-.5\wd0}}

\newcommand{\p}{\partial}
\newcommand{\epsi}{\varepsilon}
\newcommand{\dist}{\operatorname{dist}}

\begin{document}

\title[Asymptotics of harmonic functions]{Asymptotics of harmonic functions in the absence of monotonicity formulas} 

\author[Z. Li]{Zongyuan Li}
\address[Z. Li]{Department of Mathematics, Rutgers University, 110 Frelinghuysen Road, Piscataway, NJ 08854-8019, USA}
\email{zongyuan.li@rutgers.edu}
\thanks{Z. Li was partially supported by an AMS-Simons travel grant.}


\begin{abstract}
In this article, we study the asymptotics of harmonic functions. A typical method is by proving monotonicity formulas of a version of rescaled Dirichlet energy, and use it to study the renormalized solution --- the Almgren's blowup. However, such monotonicity formulas require strong smoothness assumptions on domains and operators. We are interested in the cases when monotonicity formulas are not available, including variable coefficient equations with unbounded lower order terms, Dirichlet problems on rough (non-$C^1$) domains, and Robin problems with rough Robin potentials.
\end{abstract}

\maketitle

\section{Introduction}
\begin{figure}[h] 
	\includegraphics{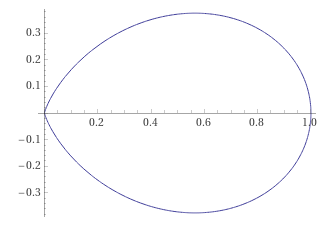} 
	\caption{Nodal set of $\mathrm{Re}(z/\log z )$}
	\label{fig-230201-0337}
\end{figure}

We discuss asymptotics of solutions to elliptic equations near both interior and boundary points. Let us start from a simple case. Consider a harmonic function $u$ in a bounded domain $\Omega \subset \bR^d$. Near an interior point $x_0 \in \Omega$, we know that $u$ is analytic:
\begin{align}
	u &= \sum_\alpha \frac{D^\alpha u(x_0)}{\alpha!} (x-x_0)^\alpha = \sum_k P_k(x-x_0) \nonumber
	\\&
	= P_N(x-x_0) + O(|x-x_0|^{N+1}). \label{eqn-230426-0353}
\end{align}
Here $P_k$ is a homogeneous harmonic polynomial of degree $k$ and $P_N$ represents the leading term. As is commonly known, expansion formulas like \eqref{eqn-230426-0353} can be useful, which are, however, not always available in the presence of variable coefficient operators or rough domains. For instance, under polar coordinates $(r,\theta)$ of $\bR^2$, consider
\begin{equation} \label{eqn-230426-0417}
	u = \mathrm{Re}\frac{r e^{i\theta}}{\log( r e^{i\theta} )}, \quad r>0, \theta \neq \pi.
\end{equation}
See Figure \ref{fig-230201-0337}. 
One can see that $u$ is harmonic in the enclosed region in Figure \ref{fig-230201-0337}, and equals to zero on the boundary given by $r =e^{\theta\tan\theta}$, except at $(r,\theta) = (1,0)$ where $u$ has a pole. Clearly it is impossible to write down an expansion like \eqref{eqn-230426-0353}, due to the log drift.

To capture the asymptotics of functions like \eqref{eqn-230426-0417}, one typically uses the ``Almgren's blowup'' --- the rescaled limit as $\lambda \rightarrow 0$ of
\begin{equation} \label{eqn-230426-0933}
	u_\lambda(\cdot) = \frac{u(\lambda\cdot)}{(\fint_{\p B_\lambda} |u|^2)^{1/2}}.
\end{equation}
For $u$ in \eqref{eqn-230426-0417}, one can simply see that $(\fint_{\p B_\lambda} |u|^2)^{1/2} \approx \lambda \log(\lambda)$ and $u_\lambda\rightarrow Cr\cos(\theta)$ as $\lambda\rightarrow 0$, where $C$ is a normalizing factor. 

Actually such convergence is guaranteed by a more general theorem --- the Almgren's monotonicity formula on convex domains. Let us describe the motivation and method. In general, one hopes to prove that the family $\{u_\lambda(\cdot)\}_{\lambda\in(0,1)}$ has one or more limits.
For this, we bound a rescaled Dirichlet energy like
\begin{equation} \label{eqn-230426-0438}
	F(r) = \frac{r D(r)}{H(r)} = \frac{r \int_{B_r} |\nabla u|^2}{\int_{\p B_r} |u|^2}.
\end{equation}
In \cite{MR574247}, Almgren observed that if $\Delta u = 0$ in $B_1$, $F(r)$ is monotonically increasing for $r\in (0,1)$. From this, $\{u_\lambda\}_{\lambda \in (0,1)}$ is uniformly bounded in $H^1$, and hence is compact in $L_2$. In literature, a quantity like \eqref{eqn-230426-0438} is usually called a (generalized) Almgren's frequencie. Its monotonicity property play an important role in blowup analysis. In this work, we are interested in three more general problems. 

\noindent\textbf{Variable coefficient equations, interior}. 
\begin{equation} \label{eqn-230426-0928}
	D_i (a_{ij}D_j u) + W_i D_i u + Vu = 0\quad \text{in} \,\, B_1,
\end{equation}
where $a_{ij}$ are symmetric, bounded, and uniformly elliptic.
In \cite{MR882069}, Garofalo-Lin proved that if $a_{ij} \in C^{0,1}, W_i, V \in L_\infty$, a modified version of $F$ in \eqref{eqn-230426-0438} is almost monotone. The condition $a_{ij} \in C^{0,1}$ cannot be improved, due to the classical counterexample in unique continuation. Later, we will discuss the cases with unbounded $W_i, V$.

\noindent\textbf{Dirichlet problem, boundary}. Suppose $\Omega \subset \bR^d$ and $0 \in \p\Omega$.
\begin{equation} \label{eqn-230426-0502}
	\begin{cases}
		\Delta u = 0\quad \text{in}\,\,\Omega\cap B_1,\\
		u = 0 \quad \text{on} \,\,\p\Omega \cap B_1.
	\end{cases}
\end{equation}
When $\Omega$ is half space or a cone, the monotonicity formula holds. For curved domains, in \cite{MR1090434, MR1363203, MR1466583}, certain variations of $F$ in \eqref{eqn-230426-0438} was proved to be almost monotone on $C^{1,1}$, convex, and $C^{1,Dini}$ domains, respectively. Some discussions on $C^1$ domains were also made in \cite{MR4544799}. It is worth mentioning that, the continuity of the normal direction $\vec{n}|_{\p\Omega}$ is essential in deriving the monotonicity formula, which is not available for rough domains, for instance general Lipschitz domains.

\noindent\textbf{Neumann and Robin problem, boundary}. Suppose $\Omega \subset \bR^d$ and $0 \in \p\Omega$.
\begin{equation}\label{eqn-230426-0531}
	\begin{cases}
		\Delta u = 0\quad \text{in}\,\,\Omega\cap B_1,\\
		\frac{\p u}{\p \vec{n}} = \eta u \quad \text{on} \,\,\p\Omega \cap B_1.
	\end{cases}
\end{equation}
Again, when $\Omega$ is half space or a cone and when $\eta = 0$ (Neumann), the monotonicity formula holds. In \cite{MR1466583, MR4104826}, this was further generalized to the case when $\p\Omega \in C^{1,1}$ and $\eta \in C^{0,1}$ (or $\eta \in W^{1,1}$ with some pointwise control on $\nabla \eta$). See also a sharp quantitative version in \cite{2021arXiv211101766L}. In all these works, the differentiability of $\eta$ cannot be dropped, which leaves the asymptotic analysis of \eqref{eqn-230426-0531} with rough $\eta$ widely open, even in the case when $\eta$ is non-negative and bounded. For instance, see the open question in \cite{MR3464049}.

\section{Alternative for motonicity formula: convergence of doubling index}

\noindent\textbf{Robin problems and variable coefficient equations.}
In a recent work, we prove the following singular set estimate.
\begin{theorem}[\cite{2023arXiv230404342L}, Theorem 1.1 (b)] \label{thm-230426-0932}
	Let $\Omega (\subset \bR^d) \in C^{1,1}$, $d\geq 2$, and $\eta \in L_p(\p\Omega)$ for some $p > 2(d-1)$. Then for any nontrivial solution $u\in H^1$ to \eqref{eqn-230426-0531}, we have
	\begin{equation*}
		\dim (\{u=|\nabla u| = 0\} \cap \Omega \cap B_{1/2}) \leq d-2.
	\end{equation*}
\end{theorem}
Such estimate relies on blowup analysis near both interior and boundary points. As mentioned before, monotonicity formulas are only proved when $\eta$ is differentiable. In \cite{2023arXiv230404342L}, we first construct an auxiliary function and reduce the problem to blowup analysis for \eqref{eqn-230426-0928} with $a_{ij} \in C^{0,1}$ and $W_i, V \in L_q$ with $q > d$. However, there is still no monotonicity formula available --- recall that the work of Garofalo-Lin \cite{MR882069} requires $W_i, V \in L_\infty$.

This requires us to design more robust methods for blowup analysis. It turns out the Federer's dimension reduction argument, which we used to prove Theorem \ref{thm-230426-0932}, only needs the following:
\begin{enumerate}
	\item a uniform $C^1$ estimate for the ``rescaled'' boundary value problems;
	\item compactness of the blowup sequence \eqref{eqn-230426-0933}, as $\lambda \rightarrow 0$;
	\item the homogeneity of the blowup limit of \eqref{eqn-230426-0933}, along subsequences.
\end{enumerate}
In \cite{2023arXiv230404342L}, (a) was achieved with the aid of the aforementioned auxiliary function and a standard $W^2_p$ regularity theory. For (b) and (c), which are typically proved via monotonicity formula, we prove the following alternative.
\begin{lemma}[\cite{2023arXiv230404342L} Lemma 4.2] \label{lem:convergence}
	Let $u \in H^1$ be a weak solution to \eqref{eqn-230426-0928} with
	\begin{equation} \label{eqn-230328-1138}
		\fint_{B_r} \left( |a_{ij} - a_{ij}(0)|^2 + r^2 |W|^2 + r^4 |V|^2 \right) \rightarrow 0,\quad \text{as}\,\,r\rightarrow 0.
	\end{equation}
	Then,
	\begin{equation*}
		\log_2 \left(\frac{\fint_{B_{2r}} |u|^2}{\fint_{B_{r}} |u|^2} \right)^{1/2} \rightarrow \bN \cup \{+\infty\}, \quad \text{as} \,\,r\rightarrow 0.
	\end{equation*}
\end{lemma}
\begin{remark}
	\begin{enumerate}
		\item The condition \eqref{eqn-230328-1138} appears naturally after scaling: if $u$ solves \eqref{eqn-230426-0928}, then $u_\lambda$ solves
		\begin{equation*}
			D_i (a_{ij}(0) D_j u_\lambda) = D_i( (a_{ij}(0) - a_{ij}(\lambda\cdot) )D_j u_\lambda) - \lambda W_i(\lambda\cdot)D_i u_\lambda - \lambda^2 V(\lambda\cdot)u_\lambda.
		\end{equation*}
		\item The condition \eqref{eqn-230328-1138} is guaranteed by $a_{ij} \in C^0$, $W_i\in L_d$, and $V \in L_{d/2}$.
		\item We will use Lemma \ref{lem:convergence} together with (SUCP) in \cite{MR1809741, MR2834777} --- if $a_{ij} \in C^{0,1}, W_i \in L_{d+\epsi}, V\in L_{d/2}$ (or $V \in L_{1+\epsi}$ when $d=2$), the limit in Lemma \ref{lem:convergence} has to be finite.
	\end{enumerate}
\end{remark}
Here in Lemma \ref{lem:convergence}, we study the doubling index
\begin{equation*}
	N(r) := \log_2 \frac{(\fint_{\p B_{2r}} |u|^2 )^{1/2}}{(\fint_{\p B_{r}} |u|^2 )^{1/2}}
\end{equation*}
instead of the frequency $F(r)$. Note that when $u$ is exactly a homogeneous polynomial of degree $k$, $N(r) \equiv k$. Hence, Lemma \ref{lem:convergence} can be interpreted as ``the existence of the limiting homogeneity''. Simple computation shows for harmonic functions, near an interior point
\begin{equation*}
	N(r) = \int_{r}^{2r} \frac{F(s)}{s} \,ds.
\end{equation*}
Hence, the monotonicity of $F$ implies the convergence of $N$, as $r\rightarrow 0$. However, the condition in Lemma \ref{lem:convergence} is much weaker than that of a monotonicity formula --- recall in \cite{MR882069}, it was required $W_i, V \in L_\infty$. Hence, we expect the conclusion of Lemma \ref{lem:convergence} can serve as a more robust tool in blowup analysis.

The proof of Lemma \ref{lem:convergence} borrows ideas of Lin-Shen \cite{MR3952693} when studying homogenization. Essentially, it is relies on fact that the monotonicity formula of harmonic functions has a rigidity property.
\begin{lemma} \label{thm-230426-1023}
	Suppose $u$ is a harmonic function in $B_1$. Then its Almgren's monotonicity function $F$ (\eqref{eqn-230426-0438}) is either strictly increasing for $r\in (0,1)$, or for some $k\in \bN$, $F\equiv k/\log 2$ and $u$ is a homogeneous harmonic polynomial of degree $k$.
\end{lemma}
From Lemma \ref{thm-230426-1023}, it can be shown that for any non-integer real number $\mu$, as $r$ decreases, after certain small scale, the doubling index $N(r)$ of $u$ can no longer jump from below $\mu$ to above. Hence, $N(r)$ is trapped near an integer, from which Lemma \ref{lem:convergence} follows.

\noindent\textbf{Dirichlet problem near a conical point.}

In a recent joint work with Dennis Kriventsov, we also study the boundary asymptotics of harmonic functions. A long-standing conjecture in boundary unique continuation asks:
\begin{conjecture}
	Suppose $u\in H^1$ is a weak solution to \eqref{eqn-230426-0502} on a Lipschitz domain $\Omega$. Then, if $\{\p u/\p\vec{n} = 0\}\cap\p\Omega$ has a postive surface measure, we must have $u \equiv 0$.
\end{conjecture}
The conjectured was proved in the case when $\Omega\in C^{1,1}, C^{1,Dini}$, and $C^1$ in \cite{MR1090434, MR1363203}, and \cite{MR1466583}, via several versions of Almgren's monotonicity formulas. For such formulas, the continuity of $\vec{n}|_{\p\Omega}$ seems inevitable, which is typically not true on general Lipschitz domains. We aim to discover the case when $\vec{n}$ is not continuous. A point $x_0 \in \p\Omega$ is called conical, if
\begin{equation*}
	\frac{(\Omega - x_0)\cap B_r}{r} \rightarrow \Gamma_{x_0} \,\, =\text{cone}.
\end{equation*}
Clearly, all differentiable $C^1$ points are conical with $\Gamma$ being the tangent plane. Moreover, any boundary point of a convex domain is conical, due to the monotonicity. In \cite{DeL}, we prove the following.
\begin{theorem}[\cite{DeL}] \label{thm-230427-0809}
	Suppose $0 \in \p\Omega$ is a conical point and $u \in H^1$ is a nontrivial solution to \eqref{eqn-230426-0502}. Then the limiting homogeneity of $u$ exists. That is,
	\begin{equation} \label{eqn-230427-0837}
				\log_2 \frac{(\fint_{\p B_{2r}\cap \Omega} |u|^2)^{1/2}}{(\fint_{\p B_{r}\cap \Omega} |u|^2)^{1/2}} \rightarrow \{\mu_j\}_j \cup \{+\infty\} \quad \text{as}\,\,r\rightarrow 0,
	\end{equation}
where $\mu_j$ are numbers determined by the spectrum of $\Delta$ on the limit cone $\Gamma$.
\end{theorem}
It is worth mentioning that, our theorem only assumes an one-point condition at $0\in\p\Omega$ --- no smoothness of $\Omega$ is needed.

\section{Uniqueness of blowup and expansion formula}
\begin{problem}
	When is the subsequence limit in \eqref{eqn-230426-0933} unique?
\end{problem}
One the one hand, naturally one may further ask.
\begin{problem}
	Does a monotonicity formula, which guarantees the existence of blowup limits, also guarantees the uniqueness of such limit?
\end{problem} 
The answer is yes when the dimension is two. This is simply due to the fact that in 2D, all eigenspaces of the Laplace operator is one-dimensional. In higher dimension, the answer is no in general. In \cite{DeL}, we constructed a convex domain $\Omega$ and a harmonic function $u$ satisfying the Dirichlet problem \eqref{eqn-230426-0502} on it. From \cite{MR1363203}, the Almgren's monotonicity formula holds due to the convexity of the domain. However, along different subsequences, the blowup limits can be different. Actually, our $\{u_\lambda\}_{\lambda\in(0,1)}$ rotates within a two-dimensional eigenspace.

On the other hand, clearly an expansion formula like \eqref{eqn-230426-0353} leads to the uniqueness of blowup limit. One can simply see that the limit has to be exactly the leading term $P_N$ upto a normalization. For Dirichlet problems, in \cite{DeL} we prove that a slightly stronger condition than ``conical'' --- H\"older conical will lead to such an expansion formula. A point $x_0 \in \p\Omega$ is called H\"older conical, if
\begin{equation*}
	\frac{\dist((\Omega-x_0)\cap B_r, \Gamma_{x_0})}{r} \leq C r^\alpha,\quad \text{for some}\,\,\alpha > 0.
\end{equation*}
\begin{theorem}[\cite{DeL}]
	Suppose $0 \in \p\Omega$ is H\"older conical. Then for any non-trivial solution $u$ to \eqref{eqn-230426-0502}, either $u = O(|x|^N)$ for all $N>0$, or there exists a $\mu_j$-homogeneous harmonic polynomial $P_{\mu_j}$ on the cone $\Gamma$, such that
	\begin{equation*}
		u(x) = P_{\mu_j}(x) + v(x),\quad \text{and} \,\, (\fint_{B_r\cap\Omega} |v|^2 )^{1/2} \leq Cr^{1+\epsi\alpha}.
	\end{equation*}
\end{theorem}
For Robin problem \eqref{eqn-230426-0531} and interior problem with variable coefficients \eqref{eqn-230426-0928}, similar property holds --- certain scaling subcritical assumptions lead to uniqueness. We expect the following:
\begin{center}
	\textit{In Lemma \ref{lem:convergence}, if we replace \eqref{eqn-230328-1138} with $a_{ij}\in C^\alpha, W_i \in L_{d+\epsi}$, and $V \in L_{d+\epsi}$, then either $u = O(|x|^N)$ for all $N>0$, or for some $k$-homogeneous harmonic polynomial $P_{k}$, we can expand $u(x) = P_k(x) + O(|x|^{k+\epsi'})$.}
\end{center}
The proof is expected to be similar to the arguments in \cite{MR1305956}. Actually, a gradient estimate is also expected for higher order terms. From these, combining the argument in \cite{MR1305956} and \cite{2023arXiv230404342L}, one can further prove the stratification of singular sets, which is stronger than our Hausdorff dimension estimate in Theorem \ref{thm-230426-0932}.

%

\end{document}